\documentclass[12pt]{amsart}
\usepackage{times,mathptm}
\usepackage{amssymb}

\renewcommand{\epsilon}{\varepsilon}
\renewcommand{\emptyset}{\varnothing}
\newcommand{\romanenumi}%
{%
}
\newcommand{\setof}[2]{\left\{\,{#1}\,:\,{#2}\:\right\}}
\newcommand{\sm}{{\setminus}}
\newcommand{\red}{\mathrm{red}}
\newcommand{\green}{\mathrm{green}}
\newcommand{\proper}[1]{\overline{#1}}

\newcommand{\toplevel}{section}

\theoremstyle{plain}
\newtheorem{Theorem}{Theorem}[\toplevel]

\newtheorem{Lemma}[Theorem]{Lemma}

\theoremstyle{definition}

\theoremstyle{remark}

\begin{document}
\title{A Suspension Lemma for Bounded Posets}
\author{J\"org Rambau}\thanks{Research at MSRI supported in part by NSF grant
  \#DMS 9022140}
\address{Konrad-Zuse-Zentrum f\"ur Informationstechnik\\
  Takustr.~7\\14195~Berlin\\
  Germany}
\email{rambau@zib.de}

\begin{abstract}
Let $P$ and $Q$ be bounded posets.
In this note, a lemma is introduced that provides a set of
sufficient conditions
for the proper part of $P$ being homotopy equivalent to the suspension
of the proper part of~$Q$.
An application of this lemma is a unified proof of the sphericity
of the higher Bruhat orders under both inclusion order (which is
a known result by \textsc{Ziegler}) and single step
inclusion order (which was not known so far).
\end{abstract}

\maketitle

\section{Introduction}

One way to draw conclusions about the homotopy type of a poset~$P$ is
to consider an order-preserving map $f$ from $P$ to another poset $Q$ the
homotopy type of which is known. If one can show that
$f$ carries a homotopy equivalence the problem is solved.
If $P$ and $Q$ are bounded then one is rather interested in the homotopy type
of the proper part $\proper{P}$ of~$P$; to take advantage of the map~$f$,
it is then usually crucial that $f: P \to Q$ restricts to
a map of the proper parts $\proper{f}: \proper{P} \to \proper{Q}$. 
However, even if this is not the case, the map $f: P \to Q$ may be exploited to
determine the homotopy type of $\proper{P}$:
in this note we present a set of sufficient conditions on $f: P \to Q$
that guarantees that $\proper{P}$ is homotopy
equivalent to the \emph{suspension} of $\proper{Q}$ (Suspension Lemma).

We apply the Suspension Lemma to show that the higher Bruhat orders by
\textsc{Manin} \& \textsc{Schechtman}~\cite{ManinSchechtman1989}
(a certain generalization of the weak Bruhat order on the symmetric group)
are spherical, no matter whether we order by \emph{inclusion} or by
\emph{single step inclusion}.

The Suspension Lemma has been applied again
in~\cite{EdelmanRambauReiner1997} to uniformly prove
the sphericity of the two (possibly different)
\emph{higher Stasheff-Tamari orders}~\cite{EdelmanReiner1996}
on the set of triangulations of a cyclic polytope.

The author would like to thank Anders Bj\"orner, Victor Reiner,
and G\"unter M.\ Ziegler for helpful discussions.

\section{The Lemma}

In this section we state and prove the Suspension Lemma.

\begin{Lemma}
Let $P,Q$ be bounded posets with $\Hat{0}_Q \neq \Hat{1}_Q$. Assume
there exist a dissection of $P$ into green elements $\green(P)$ and
red elements $\red(P)$, as well as order-preserving maps
\begin{displaymath}
  f: P \to Q \quad \text{and} \quad i,j: Q \to P
\end{displaymath}
with the following properties:
\begin{enumerate}\romanenumi
\item \label{itm:susp:Pgreenred}
  The green elements form an order ideal in~$P$.
\item \label{itm:susp:composition}
  The maps $f \circ i$ and $f \circ j$ are the identity on~$Q$.
\item \label{itm:susp:igreenjred}
  The image of $i$ is green, the image of $j$ is red.
\item \label{itm:susp:interval}
  For every $p \in P$ we have $(i \circ f)(p) \le p \le (j \circ f)(p)$.
\item \label{itm:susp:red0green1}
  The fiber $f^{-1}(\Hat{0}_Q)$ is red except for $\Hat{0}_P$, the fiber
  $f^{-1}(\Hat{1}_Q)$ is green except for~$\Hat{1}_P$.
\end{enumerate}
Then the proper part $\proper{P}$ of $P$ is homotopy equivalent to the
suspension of the proper part $\proper{Q}$ of~$Q$.
\end{Lemma}

\begin{proof}
Define
\begin{align*}
  g: 
  &\left\{
  \begin{array}{rcl}
    \proper{P} & \to     & \proper{Q \times \{ \Hat{0}, \Hat{1} \}},\\
    p       & \mapsto &
    \begin{cases}
      (f(p), \Hat{0}) & \text{if $p$ is green},\\
      (f(p), \Hat{1}) & \text{if $p$ is red};
    \end{cases}
  \end{array}
  \right.\\
\intertext{and}
  h: 
  &\left\{
  \begin{array}{rcl}
    \proper{Q \times \{ \Hat{0}, \Hat{1} \}} & \to & \proper{P},\\
    (q, \Hat{0}) & \mapsto & i(q),\\
    (q, \Hat{1}) & \mapsto & j(q).
  \end{array}
  \right.
\end{align*}

The assumptions guarantee that the above maps are well-defined and
order-preserving.
We claim that $h \circ g$
is homotopic to the identity on~$P$. In order to prove this, consider the
following carrier on the order complex $\Delta(\proper{P})$ of~$\proper{P}$.
\begin{displaymath}
  C: 
  \left\{
  \begin{array}{rcl}
    \Delta(\proper{P}) & \to & 2^{\Delta(\proper{P})},\\
    \sigma & \mapsto &
      \Delta
      \bigl(
      P_{\ge (i \circ f)(\min \sigma)}
      \cap
      P_{\le (j \circ f)(\max \sigma)}
      \cap
      \proper{P}
      \bigr).
  \end{array}
  \right.
\end{displaymath}

We claim that $C(\sigma)$ is contractible for all
$\sigma \in \Delta(\proper{P})$.
To this end, let $\sigma$ be a chain in $\proper{P}$.
If $\min \sigma$ were contained
in $f^{-1}(\Hat{0}_Q)$ and $\max \sigma$ were contained in
$f^{-1}(\Hat{1}_Q)$ then---because of
\ref{itm:susp:red0green1}---the chain $\sigma$
would have a red minimal and a green maximal element; contradiction
to~\ref{itm:susp:Pgreenred}.
Because $f \circ i$ and $f \circ j$ are the identity on $Q$,
the maps $i$ and $j$ are in particular injective. Hence, at least one of the
elements $(i \circ f)(\min \sigma)$ and
$(j \circ f)(\max \sigma)$ is contained in $\proper{P}$. 
Therefore, $C(\sigma)$ is a cone for all $\sigma \in \Delta(\proper{P})$,
thus contractible.

We further claim that the identity on $P$ and $(h \circ g)$ are both carried
by~$C$.
To see this, consider a chain $\sigma$ in $\proper{P}$ and an
element $p$ in~$\sigma$. Since
\begin{displaymath}
  (i \circ f)(\min \sigma) \stackrel{\ref{itm:susp:interval}}{\le}
  \min \sigma \le
  p \le
  \max \sigma \stackrel{\ref{itm:susp:interval}}{\le}
  (j \circ f)(\max \sigma),
\end{displaymath}
the identity on $P$ is carried by~$C$.
Because
\begin{displaymath}
  (i \circ f)(\min \sigma) \le
  (h \circ g)(\min \sigma) \le
  (h \circ g)(p) \le
  (h \circ g)(\max \sigma) \le
  (j \circ f)(\max \sigma),
\end{displaymath}
also $(h \circ g)$ is carried by~$C$.

Thus, the identity on $\proper{P}$ and $(h \circ g)$ are homotopic
by the Carrier Lemma~\cite[Lemma~10.1]{Bjoerner1995}.
Together with the fact that
$(g \circ h)$ is the identity on $Q$, this proves
that $\proper{P}$ is homotopy equivalent to
$\proper{Q \times \{ \Hat{0}, \Hat{1} \}}$.

Finally, the poset
\begin{displaymath}
  \proper{Q \times \{ \Hat{0}, \Hat{1} \}} = 
  \bigl(\proper{Q} \times \{ \Hat{0}, \Hat{1} \}\bigr) \cup
  \bigl\{ (\Hat{0}_Q, \Hat{1}), (\Hat{1}_Q, \Hat{0}) \bigr\},
\end{displaymath}
where 
\begin{displaymath}
  (\Hat{0}_Q, \Hat{1}) < \proper{Q} \times \Hat{1}
  \quad \text{and} \quad
  (\Hat{1}_Q, \Hat{0}) > \proper{Q} \times \Hat{0},
\end{displaymath}
is homeomorphic to the suspension of $\proper{Q}$ by elementary
computation rules for products and suspension of topological
spaces. Therefore, $\proper{P}$ is homotopy equivalent to
the suspension of $\proper{Q}$, as desired.
\end{proof}

\section{An Application to Higher Bruhat Orders}

In the following we present a proof for the sphericity of the
higher Bruhat orders $\mathcal{B}(n,k)$
with respect to \emph{single step inclusion order}.
Higher Bruhat orders were defined by
\textsc{Manin} \& \textsc{Schechtman}~\cite{ManinSchechtman1989}
as a generalization of the weak Bruhat order of the symmetric
group. They were further studied by \textsc{Kapranov}
\& \textsc{Voevodski}~\cite{KapranovVoevodsky1991} and
\textsc{Ziegler}~\cite{Ziegler1993}. For basic facts see these references.

For any $(k+2)$-subset $P$ of $[n]$ the set of all its $(k+1)$-subsets
is called a \emph{$(k+1)$-packet}. By abuse of notation, we denote this
$(k+1)$-packet again by~$P$.
A subset $U$ of $\tbinom{[n]}{k+1}$
is \emph{consistent} if for any $(k+1)$-packet $P$
the intersection $U \cap P$ is empty, all of $P$, or a
beginning or ending segment in the lexicographic ordering of~$P$.
For two consistent sets $U,U' \subseteq \tbinom{[n]}{k+1}$ the
\emph{single step inclusion order} is defined by
$U \le U'$ if there is a sequence $U = U_0, \dots , U_m = U'$
of consistent sets with $\#(U_i \sm U_{i-1}) = 1$ for $i = 1, \dots ,m$.

The higher Bruhat order $\mathcal{B}(n,k)$ is the set of
all consistent
subsets of $\tbinom{[n]}{k+1}$, partially ordered by single step inclusion.
In contrast to this,
$\mathcal{B}_{\subseteq}(n,k)$ is the set of all consistent subsets
of $\tbinom{[n]}{k+1}$ partially ordered by ordinary inclusion of sets
(\emph{inclusion order}).
\textsc{Ziegler}~\cite{Ziegler1993} has shown that these partial orders do
not coincide in general.

While sphericity for the inclusion order was already established
in~\cite{Ziegler1993},
the topological type of the single step inclusion order remained
an open problem. We solve this problem in the following theorem,
the proof of which works equally fine for $\mathcal{B}_{\subseteq}(n,k)$.

\begin{Theorem}
The proper part of the higher Bruhat order $\mathcal{B}(n,k)$
has the homotopy type of an $(n-k-2)$-sphere.
\end{Theorem}

\begin{proof}
We prove the theorem by induction on $n-k$.
For $n = k+1$ the higher Bruhat orders are isomorphic to
the poset $\{ \Hat{0}, \Hat{1} \}$. Therefore, $\proper{\mathcal{B}(k+1,k)}$
is the empty set, i.e., it has the homotopy type of a $(-1)$-sphere.

We show that for $n > k+1$ the conditions of the Suspension Lemma
are satisfied for 
\begin{align*}
  P &= \mathcal{B}(n,k),\\
  Q &= \mathcal{B}(n-1,k),\\
  \green(\mathcal{B}(n,k)) &= \setof{U \subseteq \tbinom{[n]}{k+1}}{
    \{ n-k, \dots , n \} \notin U },\\
  \red(\mathcal{B}(n,k)) &= \setof{U \subseteq \tbinom{[n]}{k+1}}{
    \{ n-k, \dots , n \} \in U},\\
  f&:
  \left\{
  \begin{array}{rcl}
    \mathcal{B}(n,k) & \to & \mathcal{B}(n-1,k),\\
    U & \mapsto & U \sm n := \setof{I \in U}{n \notin I};\\
  \end{array}
  \right.\\
  i&:
  \left\{
  \begin{array}{rcl}
    \mathcal{B}(n-1,k) & \to & \mathcal{B}(n,k),\\
    V & \mapsto & V;\\
  \end{array}
  \right.\\
  j&:
  \left\{
  \begin{array}{rcl}
    \mathcal{B}(n-1,k) & \to & \mathcal{B}(n,k),\\
    V & \mapsto & V \cup \setof{I \in \tbinom{[n]}{k+1}}{n \in I}.
  \end{array}
  \right.
\end{align*}

Assumptions \ref{itm:susp:Pgreenred}, \ref{itm:susp:composition},
and \ref{itm:susp:igreenjred} are obvious by the
definitions. For the inclusion order also 
\ref{itm:susp:interval} is obvious.

To prove \ref{itm:susp:interval} for the single step inclusion
order, we proceed as follows. Let $U \in \mathcal{B}(n,k)$ be a consistent
set.
We show in the sequel that $U$ can be obtained from
$(i \circ f)(U) = U \sm n$ by adding one element at a time without getting
inconsistent. Then $(i \circ f)(U) \le U$, and
we are done. (The statement about $j$ follows by taking complements.)

Let $\alpha$ an
\emph{admissible permutation} of $\tbinom{[n-1]}{k}$ corresponding
to~$U \sm n$. That is, the restriction of $\alpha$ to a $k$-packet
$P$ is the lexicographic order on $P$ if $P$ is contained in $U \sm n$; it is
the reverse lexicographic order on $P$ otherwise (see~\cite{Ziegler1993}). We
now build up $U$ from $U \sm n$ by adding the elements
$I'$ of $\setof{I \in U}{n \in I}$ in the order in which the elements
$I' \sm n$ appear in~$\alpha$. Consistency at every step
follows by construction and the fact that $\alpha$ is admissible. This
completes the proof of~\ref{itm:susp:interval}.

To see \ref{itm:susp:red0green1}, 
assume, for the sake of contradiction,
that there is a non-empty consistent set $U \in \mathcal{B}(n,k)$ with
\begin{displaymath}
  U \sm n = \emptyset 
  \quad \text{and} \quad 
  \{ n-k, \dots , n \} \notin U.
\end{displaymath}

In the following we show that every non-empty consistent set,
in particular $U$, contains at least one interval.
We call $j \in [n] \sm I$ an \emph{internal gap} of $I$ if
$\min I < j < \max I$. Note that the subsets of $[n]$ without
internal gaps are exactly the intervals.
Assume $I \in U$ has $c > 0$ internal gaps. Let $j$ be one of them.
Consider the $(k+1)$-packet $P := I \cup \{j\}$.
Since $U$ is consistent, $P \sm \min I$ or $P \sm \max I$ is
in $U$ as well. Both of them have at most $c-1$ internal gaps.
By induction we conclude that $U$ contains
at least one interval.

Since $U \sm n = \emptyset$, every element in $U$ contains~$n$.
The only interval in $\tbinom{[n]}{k+1}$ containing~$n$, however,
is $\{ n-k, \dots , n \}$. 
Hence, $\{ n-k, \dots , n \} \in U$; contradiction.

Thus,
$\Hat{0} = \emptyset$ is the only green element in
$f^{-1}(\Hat{0}) = f^{-1}(\emptyset)$.
The second statement in~\ref{itm:susp:red0green1} is again achieved
by taking complements.

Therefore, the assumptions of the Suspension Lemma are satisfied, and
$\proper{\mathcal{B}(n,k)}$ is homotopy equivalent to
the suspension of $\proper{\mathcal{B}(n-1,k)}$. This proves the theorem
by induction on~$n-k$.
\end{proof}

\providecommand{\bysame}{\leavevmode\hbox to3em{\hrulefill}\thinspace}

\end{document}